\documentclass[twocolumn]{autart}

\usepackage{savesym}
\usepackage{amsmath,amssymb,amsfonts}
\usepackage{algorithm}
\usepackage{algpseudocode}
\usepackage{autobreak}

\newtheorem{definition}{Definition}
\newtheorem{theorem}{Theorem}
\newtheorem{proposition}{Proposition}
\newtheorem{lemma}{Lemma}
\newtheorem{corollary}{Corollary}
\newtheorem{remark}{Remark}
\newtheorem{assumption}{Assumption}

\usepackage{enumerate}
\usepackage[svgnames]{xcolor}
\usepackage{tikz}
\usetikzlibrary{arrows}
\usepackage{verbatim}
\usetikzlibrary{arrows,shapes,automata}
\usetikzlibrary{fit}
\usepackage{caption}
\usepackage{subcaption}

\graphicspath{ {img/} }

\tikzset{
	treenode/.style = {align=center, inner sep=0pt, text centered,
		font=\sffamily},
	task/.style = {treenode, circle, white, font=\sffamily\bfseries, draw=black,
		fill=black, text width=1.0em},
	agent/.style = {treenode, circle, SteelBlue, draw=SteelBlue,
		text width=1.0em, thick},
}

\allowdisplaybreaks
\begin{document}
	\begin{frontmatter}
		
		\title{A Greedy and Distributable Approach to the Lexicographic Bottleneck Assignment Problem with Conditions on Exactness\thanksref{footnoteinfo}} 
		
		\thanks[footnoteinfo]{The research is funded by Defence Science and Technology Group through research agreements MyIP: 7558, MyIP: 7562 and MyIP: 9156. Corresponding author M. Khoo. Tel.:+61452523134.}
		
		\author[au1]{Mitchell Khoo}\ead{khoom1@student.unimelb.edu.au},    
		\author[au2]{Tony A. Wood}\ead{tony.wood@epfl.ch},               
		\author[au1]{Chris Manzie}\ead{manziec@unimelb.edu.au},  
		\author[au3]{Iman Shames}\ead{iman.shames@anu.edu.au}
		\address[au1]{Department of Electrical and Electronic Engineering at the University of Melbourne, Melbourne, Australia}  
		\address[au2]{Sycamore, École Polytechnique Fédérale de Lausanne (EPFL), Lausanne, Switzerland}         
		\address[au3]{CIICADA Lab, School of Engineering, Australia National University, Canberra, Australia}
		

		\begin{abstract}                          
			Solving the Lexicographic Bottleneck Assignment Problem (LexBAP) typically relies on centralised computation with order $\mathcal{O}(n^4)$ complexity. We consider the Sequential Bottleneck Assignment Problem (SeqBAP), which yields a greedy solution to the LexBAP and discuss the relationship between the SeqBAP, the LexBAP, and the Bottleneck Assignment Problem (BAP). In particular, we reexamine tools used to analyse the structure of the BAP, and apply them to derive an $\mathcal{O}(n^3)$ algorithm that solves the SeqBAP. We show that the set of solutions of the LexBAP is a subset of the solutions of the SeqBAP and analyse the conditions for which the solutions sets are identical. Furthermore, we provide a method to verify the satisfaction of these conditions. In cases where the conditions are satisfied, the proposed algorithm for solving the SeqBAP solves the LexBAP with computation that has lower complexity and can be distributed over a network of computing agents. The applicability of the approach is demonstrated with a case study where mobile robots are assigned to goal locations.
		\end{abstract}
		
	\end{frontmatter}

\section{Introduction}
In multi-agent systems, assignment problems arise when a set of tasks must be allocated to a set of agents, where each allocation of a task to an agent incurs a cost. Reviews on assignment problems with different objectives are found in \cite{assignprob,survey,taxonomy}. The Bottleneck Assignment Problem (BAP) is an assignment problem with the objective of minimising the costliest allocation. It appears in time-critical problems, e.g., when agents carry out tasks simultaneously and must complete all tasks in minimum time. For instance, the BAP arises in~\cite{decoys}, where the goal is to minimise the worst-case positioning time of decoys.

A centralised algorithm relies on a single decision-maker to aggregate all information from agents to compute the solution. In contrast, a distributed algorithm is one where computation is distributed over agents. Several centralised algorithms to solve the BAP have been proposed in \cite{BAP1,BAP2,BAP3,BAP4}, while a distributed algorithm to solve the BAP is presented in \cite{distBAP}.

Two special cases of the BAP are the Lexicographic Bottleneck Assignment Problem (LexBAP) \cite{lex1} and the Sequential Bottleneck Assignment Problem (SeqBAP) \cite{tony,seq1}. The former focuses not only on minimising the costliest allocation, but also the second costliest, and the third costliest, etc. The latter, which is further elaborated in this paper, is a greedy reformulation of the former, where each allocation is chosen sequentially with no regard to the effect on later choices. It is common to reformulate greedy versions of assignment problems, e.g., a distributed greedy algorithm for an assignment problem with submodular utility functions is presented in~\cite{distgreed}, while a greedy rescheduling of an exact LexBAP solution given runtime uncertainty is proposed in~\cite{seq2}.

Although a solution of the SeqBAP may not be a solution of the LexBAP in general, there are conditions for which the solution sets of the SeqBAP and LexBAP coincide. This motivates the derivation of an algorithm that solves the SeqBAP. Existing algorithms for solving the LexBAP, e.g., in \cite{assignprob,lex1}, are not amenable to a distributed implementation. In this paper, we develop an approach to solve the SeqBAP that can be implemented with distributed computation and provides certificates when the resulting assignment is also a solution to the LexBAP.

This paper extends the work in~\cite{distBAP,struct,tony} in the following ways. The distributed algorithm for finding a solution to the BAP in~\cite{distBAP} is extended to develop an algorithm to produce a solution to the SeqBAP that is similarly amenable to a distributed implementation. Applying a BAP algorithm ``off-the-shelf'' to solve the SeqBAP was first proposed in~\cite{tony,seq1}. However, we henceforth refer to such an approach as a naive one as it does not exploit structure in the SeqBAP.

In~\cite{struct}, tools are introduced to identify and exploit structure of the BAP. In particular, two BAPs with two distinct subsets of agents and tasks are solved separately and the solutions to the separate problems are used to solve the BAP over the combined sets of agents and tasks efficiently. We draw a parallel to~\cite{struct} by similarly analysing the structure of the SeqBAP, and exploit it to modify the aforementioned naive SeqBAP approach so that the resulting novel approach has lower worst-case complexity. We prove that the novel approach has lower theoretical worst-case complexity than the naive approach and demonstrate that it also has lower empirical complexity in a case study.

Conditions for which a solution of the LexBAP and SeqBAP are unique and equal were introduced in~\cite{tony}. We extend this analysis further and determine the full relationship between the SeqBAP and the LexBAP.  We show that the solution set of the LexBAP is a subset of the solutions of the SeqBAP, which in turn is a subset of the solution set of the BAP.

By establishing that the SeqBAP is a greedy approximation of the LexBAP with certifiable conditions for exactness, the main contribution of this paper is the derivation an efficient algorithm to solve the LexBAP. The benefits of the proposed algorithm compared to existing literature are two-fold. Firstly, solving the LexBAP according to~\cite{assignprob} has a worst-case complexity of $\mathcal{O}(n^4)$, where $n$ is number of allocations of tasks to agents that have to be established. By exploiting structure of the SeqBAP, we instead present an algorithm for the SeqBAP that has a worst-case complexity of $\mathcal{O}(n^3)$. Secondly, this algorithm is intrinsically distributable and can be implemented with computation distributed over agents.

The rest of the paper is organised as follows. In Section~\ref{sec:theory}, background graph theory is introduced. In Section~\ref{sec:problems}, the BAP, LexBAP, and the SeqBAP are formulated. In Section~\ref{sec:structs}, the relationship between these assignment problems is analysed and tools to exploit their structure are derived. In Section~\ref{sec:alg}, the distributable algorithm for solving the SeqBAP is presented. In Section~\ref{sec:case}, a numerical case study is provided to demonstrate the implementation of the algorithm.

\section{Graph Theoretical Definitions} \label{sec:theory}


Consider an undirected graph $\mathcal{G}=(\mathcal{V,E})$, where $\mathcal{V}$ is a set of vertices and $\mathcal{E}$ is a set of edges. Given $\mathcal{G}$, we consider the following definitions, as in \cite{distBAP,hk_alg}.

\begin{definition}[Maximum Cardinality Matching] \label{def:matching}
	A matching $\mathcal{M}$ of $\mathcal{G}$ is a subset of edges $\mathcal{M}\subseteq \mathcal{E}$ such that no vertex $v\in\mathcal{V}$ is incident to more than one edge in $\mathcal{M}$. A Maximum Cardinality Matching (MCM) $\mathcal{M}_{max}$ is a matching of $\mathcal{G}$ with maximum cardinality.
\end{definition}

\begin{definition}[Neighbours] \label{def:neigh}
	The set of neighbours of vertex $v\in\mathcal{V}$ in $\mathcal{G}$ is $N(\mathcal{G},v):=\{k|\{v,k\}\in \mathcal{E}\}$.
\end{definition}

\begin{definition}[Path] \label{def:path}
	Let distinct vertices $v_1, v_2, ..., v_{l+1}\in \mathcal{V}$ be such that for $k\in\{1,2,...,l\}$, $v_{k+1}\in N(\mathcal{G},v_k)$. The set of edges $\mathcal{P}=\{ \{v_k,v_{k+1}\}|k\in\{1,2,...,l\}\}$ is a path between $v_1$ and $v_{l+1}$, with length $l$.
\end{definition}

\begin{definition}[Diameter]
	Let $h_{ij}$ be the length of the shortest path between vertices $i,j\in\mathcal{V}$. The diameter $D$ of graph $\mathcal{G}$ is $D:=\max_{i,j\in\mathcal{V}}h_{ij}$.
\end{definition}

\begin{definition}[Alternating path] \label{def:alt}
	Given a matching $\mathcal{M}$ and a path $\mathcal{P}$, $\mathcal{P}$ is an alternating path relative to $\mathcal{M}$ if and only if each vertex that is incident to an edge in $\mathcal{P}$ is incident to no more than one edge in $\mathcal{P}\cap\mathcal{M}$ and no more than one edge in $\mathcal{P}\backslash\mathcal{M}$.
\end{definition}

A path is a set of edges and its elements are unordered. However, if the elements of an alternating path relative to $\mathcal{M}$ were arranged in a sequence $\{v_1,v_2\}$, $\{v_2,v_3\}$, ..., $\{v_{l-1},v_l\}$, $\{v_l,v_{l+1}\}$, then the edges in the sequence alternate between edges in $\mathcal{M}$ and edges not in $\mathcal{M}$.

\begin{definition}[Free vertex] \label{def:free}
	Given a matching $\mathcal{M}$, a vertex $v\in \mathcal{V}$ is free if and only if for all $w\in \mathcal{V}$, $\{v,w\}\notin \mathcal{M}$.
\end{definition}

\begin{definition}[Augmenting path] \label{def:aug}
	Given a matching $\mathcal{M}$ and a path $\mathcal{P}$ between vertices $v_1$ and $v_{l+1}$, $\mathcal{P}$ is an augmenting path relative to $\mathcal{M}$ if and only if $\mathcal{P}$ is an alternating path relative to $\mathcal{M}$ and $v_1$ and $v_{l+1}$ are both free vertices.
\end{definition}


%


\section{Assignment Problem Formulations} \label{sec:problems}
Consider a bipartite graph $\mathcal{G}_b=(\mathcal{V}_b,\mathcal{E}_b)$, with vertex set $\mathcal{V}_b=\mathcal{A}_b\cup \mathcal{B}_b$ and edge set $\mathcal{E}_b \subseteq \{\{i,j\}|i\in\mathcal{A}_b,j\in\mathcal{B}_b\}$, where $\mathcal{A}_b$ is a set of agents and $\mathcal{B}_b$ is a set of tasks such that $\mathcal{A}_b\cap \mathcal{B}_b=\emptyset$. Let $\mathcal{C}(\mathcal{G}_b)$ be the set of all MCMs of $\mathcal{G}_b$. Let the function $w:\mathcal{E}_b\rightarrow \mathbb{R}$ map the edges of $\mathcal{G}_b$ to real-valued weights.

Consider an MCM $\mathcal{M}\in \mathcal{C}(\mathcal{G}_b)$. Let $n=|\mathcal{M}|$. Without loss of generality, we denote the edges in $\mathcal{M}$ as $\mathcal{M}=\{\{i_1,j_1\},\{i_2,j_2\}, ..., \{i_n,j_n\}\}$ and assume $w(\{i_1,j_1\})$ $\geq$ $w(\{i_2,j_2\})$ $\geq$ $...$ $\geq$ $w(\{i_n,j_n\})$. We define the ordered tuple of weights $\mathcal{T}(\mathcal{M}):=(w(\{i_1,j_1\}), w(\{i_2,j_2\}), ..., w(\{i_n,j_n\})$ and the $k$th element of the tuple $\mathcal{T}_k(\mathcal{M}):=w(\{i_k,j_k\})$, for $k\in\{1,2,...,n\}$.


\subsection{The Bottleneck Assignment Problem}

The BAP for graph $\mathcal{G}_b$ is formulated as
\begin{subequations} \label{eq:BAP}
\begin{align}
\nonumber \text{BA}&\text{P}:\\
&\text{Find }\mathcal{M}\in\mathcal{C}(\mathcal{G}_b), \text{ s.t. }\\
\nonumber&\forall \mathcal{M}'\in\mathcal{C}(\mathcal{G}_b)\backslash\{\mathcal{M}\},\\
&\mathcal{T}_1(\mathcal{M})\leq \mathcal{T}_1(\mathcal{M}').
\end{align}
\end{subequations}
We define the bottleneck weight, bottleneck assignment and the bottleneck edge of a bipartite graph $\mathcal{G}_b$.

\begin{definition}[Bottleneck weight] \label{def:botweight}
	The bottleneck weight of a bipartite graph $\mathcal{G}_b$ is defined as $\mathcal{W}(\mathcal{G}_b):=\min_{\mathcal{M}\in \mathcal{C}(\mathcal{G}_b)}\max_{\{i,j\}\in \mathcal{M}} w(\{i,j\})$.
\end{definition}

\begin{definition}[Bottleneck assignment]\label{def:botass}
	The set of bottleneck assignments of $\mathcal{G}_b$ is $\mathcal{S}(\mathcal{G}_b):= \{\mathcal{M}\in\mathcal{C}(\mathcal{G}_b)|\mathcal{T}_1(\mathcal{M})=\mathcal{W}(\mathcal{G}_b)\}$.
\end{definition}


\begin{definition}[Bottleneck edge] \label{def:botedge}
	Given any bottleneck assignment $\mathcal{M}\in\mathcal{S}(\mathcal{G}_b)$, any edge\linebreak $e\in \arg\max_{\{i,j\} \in\mathcal{M}} w(\{i,j\})$ is a bottleneck edge of $\mathcal{G}_b$.
\end{definition}

\begin{definition}[Matching-sublevel set] \label{def:sublevel}
	Given a bipartite graph $\mathcal{G}_b=(\mathcal{V}_b,\mathcal{E}_b)$ with a matching $\mathcal{M}$, define a matching-sublevel set $\psi(\mathcal{G}_b,\mathcal{M}):=\{e\in\mathcal{E}_b|w(e)\leq \max_{e'\in\mathcal{M}}w(e')\}$ and a strict matching-sublevel set $\psi_S(\mathcal{G}_b,\mathcal{M}):=\{e\in\mathcal{E}_b| w(e)< \max_{e'\in\mathcal{M}}w(e')\}$.
\end{definition}

The following definition was first introduced in~\cite{struct} to capture the structure in the BAP.

\begin{definition}[Critical bottleneck edge] \label{def:crit}
	Let $\mathcal{M}$ be an MCM of graph $\mathcal{G}_b$. Edge $e_c$ is a critical bottleneck edge of $\mathcal{G}_b$ relative to $\mathcal{M}$ if and only if $e_c\in\arg\max_{e\in\mathcal{M}} w(e)$ and $\phi(\mathcal{G}_b,\mathcal{M})\backslash\{e_c\}$ does not contain an augmenting path relative to $\mathcal{M}\backslash \{e_c\}$, where $\phi(\mathcal{G}_b,\mathcal{M}):=\mathcal{M}\cup\psi_S(\mathcal{G}_b,\mathcal{M})$ and $\psi_S(\cdot)$ is defined in Definition~\ref{def:sublevel}.
\end{definition}

In~\cite{struct}, it is shown that a critical bottleneck edge is also a bottleneck edge. A critical bottleneck edge is used to identify if an MCM is a solution to~(\ref{eq:BAP}). From Definition~\ref{def:crit}, we observe that the removal of a critical bottleneck edge results in the non-existence of an augmenting path. Applied together with Berge's Theorem~\cite{berge}, it is concluded that finding a critical bottleneck edge in an MCM $\mathcal{M}$ implies that $\mathcal{M}$ is a solution to~(\ref{eq:BAP}).

\subsection{The Lexicographic Bottleneck Assignment Problem}
A special case of the BAP is the LexBAP. The LexBAP is formulated as
\begin{subequations} \label{eq:LBAP}
	\begin{align}
	\text{Le}&\text{xBAP}:\nonumber \\
	&\text{Find }\mathcal{M}\in\mathcal{C}(\mathcal{G}_b), \text{ s.t.}\\
	\nonumber&\forall \mathcal{M}'\in\{\tilde{\mathcal{M}}\in\mathcal{C}(\mathcal{G}_b)|\mathcal{T}(\tilde{\mathcal{M}})\neq \mathcal{T}(\mathcal{M})\},\\
	\nonumber&\exists k\in\{z\in\mathbb{Z}^+|z\leq n\}, \mathcal{T}_k(\mathcal{M})<\mathcal{T}_k(\mathcal{M}'),\\
	&\forall l\in\{z\in\mathbb{Z}^+|z<k\}, \mathcal{T}_l(\mathcal{M})=\mathcal{T}_l(\mathcal{M}'),
	\end{align}
\end{subequations}
where $\mathbb{Z}^+$ is the set of all strictly positive integers. For example, the tuple $\mathcal{T}(\mathcal{M}_1)=(5,3,3,3)$ is lexicographically smaller than $\mathcal{T}(\mathcal{M}_2)=(5,4,3,2)$ because its second element is smaller and their first elements have equal weight. The objective of the LexBAP is to choose a lexicographically minimal MCM. Note that while the largest edge weights of $\mathcal{M}_1$ and $\mathcal{M}_2$ are equal to 5, they may not correspond to the same edge. In the greedy approach below, this ambiguity could lead to $\mathcal{M}_2$ being produced as a solution instead of $\mathcal{M}_1$.

\subsection{A Greedy Solution to the Lexicographic Bottleneck Assignment Problem}

A greedy approach to the LexBAP involves sequentially solving $n$ BAPs, i.e., sequentially choosing the edges corresponding to $\mathcal{T}_1(\mathcal{M})$, $\mathcal{T}_2(\mathcal{M})$, ..., $\mathcal{T}_n(\mathcal{M})$ one at a time. This does not in general produce a solution to the LexBAP. Each time an edge is selected, it affects the remaining selections in the sequence. The greedy selection of edges may involve an arbitrary choice between edges with the same weight, e.g., $\mathcal{M}_1$ and $\mathcal{M}_2$ from the example in the previous section that both have an edge with weight 5. Therefore, the greedy approach can potentially produce a suboptimal solution to the LexBAP.


A particular greedy approximation of the LexBAP is introduced in \cite{tony}, which we henceforth refer to as SeqBAP. To formulate the SeqBAP we first introduce the notion of the price of absence of an edge. To this end, let $\mathcal{F}(\mathcal{G}_b)$ be the cardinality of an MCM of $\mathcal{G}_b$, i.e., given any MCM $\mathcal{M}\in\mathcal{C}(\mathcal{G}_b)$, $\mathcal{F}(\mathcal{G}_b)=|\mathcal{M}|$.

\begin{definition}[Price of Absence]
	Given a bipartite graph $\mathcal{G}_b=(\mathcal{V}_b,\mathcal{E}_b)$, and an edge $e\in\mathcal{E}_b$ such that $\mathcal{F}((\mathcal{V}_b,\mathcal{E}_b))=\mathcal{F}((\mathcal{V}_b,\mathcal{E}_b\backslash\{e\}))$, the price of the absence of an edge $e$ is \begin{equation*} P(\mathcal{G}_b,e):=\mathcal{W}((\mathcal{V}_b,\mathcal{E}_b\backslash\{e\}))-\mathcal{W}((\mathcal{V}_b,\mathcal{E}_b)),\end{equation*} where $\mathcal{W}(\cdot)$ is defined in Definition~\ref{def:botweight}. If the removal of an edge $e\in\mathcal{E}_b$ changes the cardinality of an MCM, i.e., $\mathcal{F}((\mathcal{V}_b,\mathcal{E}_b))\neq \mathcal{F}((\mathcal{V}_b,\mathcal{E}_b\backslash\{e\}))$, then the price of absence is defined to be $+\infty$.
\end{definition}
The price of absence is always non-negative. It is a measure of how much the bottleneck weight of a graph increases in the absence of a given edge. The so-called robustness margin in~\cite{tony} is a special case of the price of absence. Its value can be used to quantify the sensitivity of an assignment solution to perturbations of edge weights as first studied in~\cite{LAP2}. This type of sensitivity information is exploited in~\cite{tony} to guarantee collision avoidance of mobile agents that are assigned to different destinations. We show that the price of absence is also related to a critical bottleneck edge in the following proposition.


\begin{proposition} \label{prop:POA}
	Given a bipartite graph $\mathcal{G}_b=(\mathcal{V}_b,\mathcal{E}_b)$, let the MCM $\mathcal{M}\in\mathcal{S}(\mathcal{G}_b)$ be a bottleneck assignment of $\mathcal{G}_b$ and let edge $e_b\in\mathcal{M}$ be a bottleneck edge of $\mathcal{G}_b$. If $e_b$ has a positive price of absence, then it is a critical bottleneck edge of $\mathcal{G}_b$ relative to $\mathcal{M}$.
\end{proposition}
\begin{pf}
	Assume bottleneck edge $e_b$ is not a critical bottleneck edge of $\mathcal{G}_b$ relative to $\mathcal{M}$. By Definition~\ref{def:crit}, $\phi(\mathcal{G}_b,\mathcal{M})\backslash \{e_b\}$ contains an augmenting path $\mathcal{P}$ relative to $\mathcal{M}\backslash \{e_b\}$, where $\phi(\mathcal{G}_b,\mathcal{M}):=\mathcal{M}\cup\psi_S(\mathcal{G}_b,\mathcal{M})$ and $\psi_S(\cdot)$ is given in Definition~\ref{def:sublevel}. By Berge's Theorem~\cite{berge}, $\mathcal{M}':=\mathcal{M}\backslash \{e_b\}\oplus \mathcal{P}$ is an MCM of $\mathcal{G}_b$, where the operator $\oplus$ denotes symmetric difference. The weight of all edges in $\mathcal{M}'$ must be smaller than or equal to $w(e_b)$ because $\mathcal{M}'\subseteq \phi(\mathcal{G}_b,\mathcal{M})\backslash \{e_b\}$. This holds as $\mathcal{M}'$ is formed from the symmetric difference of $\mathcal{M}\backslash \{e_b\}\subseteq \phi(\mathcal{G}_b,\mathcal{M})\backslash \{e_b\}$ and $\mathcal{P}\subseteq \phi(\mathcal{G}_b,\mathcal{M})\backslash \{e_b\}$. However, $\phi(\mathcal{G}_b,\mathcal{M})\backslash \{e_b\}\subset \mathcal{E}_b\backslash \{e_b\}$ implies that $\mathcal{E}_b\backslash \{e_b\}$ must also contain the MCM $\mathcal{M}'$, i.e., $\mathcal{W}((\mathcal{V}_b,\mathcal{E}_b\backslash \{e_b\}))\ngtr \mathcal{W}((\mathcal{V}_b,\mathcal{E}_b))$. $\hfill\qed$
\end{pf}

The SeqBAP is constructed by sequentially choosing the bottleneck edge with maximal price of absence and removing the corresponding bottleneck agent and task from the graph. Thus, it is a greedy solution to LexBAP. Given a set of weighted edges $\mathcal{E}\subseteq\mathcal{E}_b$, we define
\begin{equation} \label{eq:mathl}
\mathcal{L}(\mathcal{E}):=\arg\max_{\{i,j\}\in\mathcal{E}} w(\{i,j\}).
\end{equation}%
The SeqBAP is formulated as
\begin{subequations} \label{eq:SBAP}
	\begin{align}
	\text{Se}&\textrm{qBAP}:\nonumber\\
	\label{eq:seqedge}&\text{Find } \{\{i_1,j_1\},\{i_2,j_2\},...,\{i_n,j_n\}\}\in\mathcal{C}(\mathcal{G}_b),\\
	&\text{s.t. } \forall k\in\{z\in\mathbb{Z}^+|z\leq n\},\nonumber\\
	\label{eq:seqedgek}&\{i_k,j_k\}\in\arg\!\max_{\{i,j\}\in \mathcal{L}(\mathcal{M}^k)}P(\mathcal{G}^k,\{i,j\}),\\
	\label{eq:M_k}&\mathcal{M}^k \in \arg\!\min_{\mathcal{M}\in \mathcal{C}(\mathcal{G}^k)}\max_{\{i,j\}\in \mathcal{M}} w(\{i,j\}),\\ 
	\label{eq:G_k}&\mathcal{G}^k = (\mathcal{A}^k\cup\mathcal{B}^k,\mathcal{E}^k),\\
	\label{eq:E_k}&\mathcal{E}^k =\{\{i,j\}\in \mathcal{E}_b|i\in\mathcal{A}^k,j\in\mathcal{B}^k\},\\
	&\text{where } \forall k\in\{z\in\mathbb{Z}^+|2\leq z\leq n\}, \nonumber\\
	\label{eq:A_1}&\mathcal{A}^1 = \mathcal{A}_b, \mathcal{B}^1=\mathcal{B}_b,\\
	\label{eq:A_k}&\mathcal{A}^k = \mathcal{A}^{k-1}\backslash \{i_{k-1}\},\\
	\label{eq:B_k}&\mathcal{B}^k = \mathcal{B}^{k-1}\backslash \{j_{k-1}\}.
	\end{align}
\end{subequations}

\section{Structure of the BAP, the LexBAP, and the SeqBAP} \label{sec:structs}

We exploit the structure of the BAP, the LexBAP, and the SeqBAP in two ways. In order to derive an efficient method to solve the SeqBAP, we consider the role of all edges with positive price of absence within the BAP, the LexBAP, and the SeqBAP, and exploit this to circumvent the need to find the edge with maximum price of absence as indicated in~(\ref{eq:seqedgek}). To guarantee that a solution to the SeqBAP is an exact solution to the LexBAP, we consider the relationship between their solution sets. In particular, the subsections are organised as follows.

In Section~\ref{subsec:identify}, we apply the tools introduced in Section~\ref{sec:problems} to efficiently identify edges with positive price of absence. In Section~\ref{subsec:solutions}, we show that finding a SeqBAP solution does not require computation of the explicit price of absence for each edge in~(\ref{eq:seqedgek}) and only requires identification of edges with positive price. Additionally, in Section~\ref{subsec:solutions} we analyse conditions for a solution to SeqBAP to be an exact solution to the LexBAP. These results lead to the derivation of a SeqBAP algorithm in Section~\ref{sec:alg}, which serves as a greedy approach to finding a solution to the LexBAP with exactness guarantees.


\subsection{Identifying Edges with Positive Price of Absence} \label{subsec:identify}

The following proposition establishes that an edge with positive price of absence appears in all solutions to the BAP. This result is a generalisation of a property of robustness margins proven in~\cite{tony}.


\begin{proposition} \label{prop:priceBOT}
	Consider a bipartite graph $\mathcal{G}_b=(\mathcal{V}_b,\mathcal{E}_b)$. If an edge $e_p\in\mathcal{E}_b$ has a positive price of absence, then $e_p$ is an element of every bottleneck assignment of $\mathcal{G}_b$, i.e., for all $\mathcal{M}\in\mathcal{S}(\mathcal{G}_b), e_p\in\mathcal{M}$.
\end{proposition}
\begin{pf}
	Consider an arbitrary edge $e_p\in\mathcal{E}_b$ with $P(\mathcal{G}_b,e_p)=\mathcal{W}((\mathcal{V}_b,\mathcal{E}_b\backslash \{e_p\}))-\mathcal{W}((\mathcal{V}_b,\mathcal{E}_b))>0$. Assume for the sake of contradiction that there exists an MCM $\mathcal{M}\in\mathcal{S}(\mathcal{G}_b)$ such that $e_p\notin \mathcal{M}$. However, this implies that $\mathcal{W}((\mathcal{V}_b,\mathcal{E}_b\backslash \{e_p\}))=\mathcal{W}((\mathcal{V}_b,\mathcal{E}_b))$, which contradicts the assumption that $P(\mathcal{G}_b,e_p)>0$. $\hfill\qed$
\end{pf}

The following corollary follows from Proposition~\ref{prop:priceBOT} for a set of edges with positive price of absence.

\begin{corollary} \label{cor:E}
	Given a bipartite graph $\mathcal{G}_b=(\mathcal{V}_b,\mathcal{E}_b)$, let MCM $\mathcal{M}\in\mathcal{S}(\mathcal{G}_b)$ be a bottleneck assignment of $\mathcal{G}_b$ and let $E' := \{e\in\mathcal{M} | P(\mathcal{G}_b,e)>0\}$ be the set of edges with positive price of absence. The set $E'$ is a subset of every bottleneck assignment of $\mathcal{G}_b$, i.e., for all $\mathcal{M}'\in\mathcal{S}(\mathcal{G}_b), E'\subseteq\mathcal{M}'$.
\end{corollary}

The following theorem provides a property of an edge $e_p$ that can be exploited to identify whether $e_p$ has positive price of absence.

\begin{theorem} \label{theorem:pos}
	Given a bipartite graph $\mathcal{G}_b=(\mathcal{V}_b,\mathcal{E}_b)$, let MCM $\mathcal{M}\in\mathcal{S}(\mathcal{G}_b)$ be a bottleneck assignment of $\mathcal{G}_b$. The edge $e_p\in\mathcal{E}_b$ has positive price of absence if and only if there does not exist an augmenting path in $\psi(\mathcal{G}_b,\mathcal{M})\backslash\{e_p\}$ relative to $\mathcal{M}\backslash\{e_p\}$, where the matching-sublevel set $\psi(\cdot)$ is defined in Definition~\ref{def:sublevel}.
\end{theorem}
\begin{pf}
	First, we prove necessity. Assume $e_p$ has positive price of absence. Assume for the sake of contradiction that there exists an augmenting path $\mathcal{P}'\subseteq \psi(\mathcal{G}_b,\mathcal{M})\backslash\{e_p\}$ relative to $\mathcal{M}\backslash\{e_p\}$. By Berge's Theorem~\cite{berge}, $\mathcal{M}'=\mathcal{P}'\oplus \mathcal{M}\backslash \{e_p\}$ is another MCM of $\mathcal{G}_b$, where the operator $\oplus$ denotes symmetric difference. However, $\mathcal{M}'$ does not contain $e_p$ because it is the symmetric difference of two sets of edges that both do not contain $e_p$. This contradicts Proposition~\ref{prop:priceBOT} because $e_p$ is an edge with positive price of absence and all solutions to~(\ref{eq:BAP}) must contain $e_p$.
	
	Next, we prove sufficiency. Assume there does not exist an augmenting path $\mathcal{P}'\subseteq \psi(\mathcal{G}_b,\mathcal{M})\backslash\{e_p\}$ relative to $\mathcal{M}\backslash\{e_p\}$. Then, it is not possible to construct an MCM of $\mathcal{G}_b$ from $\psi(\mathcal{G}_b,\mathcal{M})\backslash\{e_p\}$. In other words, without edge $e_p$, an MCM of $\mathcal{G}_b$ must contain at least one edge with weight larger than the weights of the edges in $\mathcal{M}$. Thus, $e_p$ has positive price of absence. $\hfill\qed$
\end{pf}

The following corollary follows from Theorem~\ref{theorem:pos} and further illustrates the property of an edge with positive price of absence and its role within the structure of a matching-sublevel set.

\begin{corollary} \label{prop:price}
	Given a bipartite graph $\mathcal{G}_b=(\mathcal{V}_b,\mathcal{E}_b)$, let MCM $\mathcal{M}\in\mathcal{S}(\mathcal{G}_b)$ be a bottleneck assignment of $\mathcal{G}_b$. The edge $e_p=\{a_p,b_p\}\in\mathcal{E}_b$ has positive price of absence if and only if the path $\mathcal{P}=\{e_p\}$ is a unique alternating path in the matching-sublevel set $\psi(\mathcal{G}_b,\mathcal{M})$ relative to $\mathcal{M}$ between $a_p$ and $b_p$.
\end{corollary}

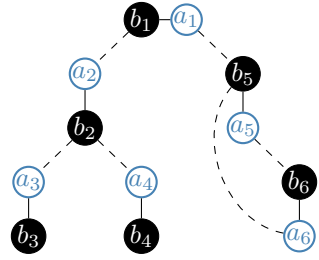
\begin{figure}[thpb]
	\centering
	\begin{tikzpicture}[-,>=stealth',level/.style={sibling distance = 2.5cm,
		level distance = 1.2cm},scale=0.6]
	\node [task] at (0,0) (b1) {$b_1$}
	child{  node [agent] (a2) {$a_2$} edge from parent[dashed]
		child[solid]{ node [task] (b2) {$b_2$} edge from parent[solid]
			child{  node [agent] (a3) {$a_3$} edge from parent[dashed]
				child[solid]{ node [task] (b3) {$b_3$} edge from parent[solid] node[left] {}
				} node[left] {}
			} 
			child{ node [agent] (a4) {$a_4$} edge from parent[dashed]
				child[solid]{ node [task] (b4) {$b_4$} edge from parent[solid] node[left] {}
				} node[right] {}
			} node[left] {}
		} node[left] {}
	}
	child{ node [] {} edge from parent[draw=none]
	}
	;
	
	\node [agent] at (1,0) (a1) {$a_1$}
	child{ node [] {} edge from parent[draw=none]
	}
	child{  node [task] (b5) {$b_5$} edge from parent[dashed]
		child[solid]{ node [agent] (a5) {$a_5$} edge from parent[solid]
			child{ node [] {} edge from parent[draw=none]
			}
			child{  node [task] (b6) {$b_6$} edge from parent[dashed]
				child[solid]{ node [agent] (a6) {$a_6$} edge from parent[solid] node[right] {}
				} node[right] {}
			} node[right] {}
		} node[right] {}
	}
	;
	\path [draw] (b1) -- node [label=below:{}] {} (a1);
	\path [draw,dashed] (a6) edge[bend left=60] node[right] {} (b5);
	\end{tikzpicture}
	\caption{Illustration of Theorem~\ref{theorem:pos}: The lines in the figure represent edges in the matching-sublevel set $\psi(\mathcal{G}_b,\mathcal{M})$ for a graph $\mathcal{G}_b$ and an MCM $\mathcal{M}$. White nodes represent agents $\mathcal{A}_b$ and black nodes represent tasks $\mathcal{B}_b$. The edges in the MCM $\mathcal{M}$ are shown as solid lines and the edges not in $\mathcal{M}$ are shown as dashed lines. In this example, all edges in $\mathcal{M}$ have positive price of absence except for $\{a_5,b_5\}$ and $\{a_6,b_6\}$. We consider edge $\{a_2,b_2\}$ with positive price of absence. There exists only one alternating path between $a_2$ and $b_2$, i.e., the path $\mathcal{P}=\{\{a_2,b_2\}\}$. If this edge were removed from the matching-sublevel set, there would not exist an augmenting path between $a_2$ and $b_2$. No MCM of $\mathcal{G}_b$ can be constructed using edges contained in $\psi(\mathcal{G}_b,\mathcal{M})\backslash\{\{a_2,b_2\}\}$. In contrast if edge $\{a_5,b_5\}$ were removed, there exists an augmenting path $\mathcal{P}=\{\{a_6,b_5\},\{a_6,b_6\},\{a_5,b_6\}\}$ between $a_5$ and $b_5$ relative to $\mathcal{M}\backslash\{\{a_5,b_5\}\}$. The matching $\mathcal{M}'=(\mathcal{M}\backslash\{\{a_5,b_5\}, \{a_6,b_6\}\}) \cup \{\{a_6,b_5\}, \{a_5,b_6\}\}$ is an MCM of $\mathcal{G}_b$ with all edges in $\mathcal{M}'$ being elements of the set $\psi(\mathcal{G}_b,\mathcal{M})\backslash\{\{a_5,b_5\}\}$.} \label{fig:tree2}
\end{figure}

Fig.~\ref{fig:tree2} illustrates the results from Theorem~\ref{theorem:pos} and Corollary~\ref{prop:price} within a simple example.

Selecting edges for an MCM that solves the SeqBAP in~(\ref{eq:seqedge}) only requires checking strict positivity of the price of absence, not the explicit value. In order to evaluate the price of absence of an edge explicitly, an additional BAP would have to be solved. Theorem~\ref{theorem:pos} provides a method to determine if an edge has positive price of absence without evaluating this extra BAP by instead searching for an augmenting path. Finding one augmenting path is less complex than solving one BAP, which itself requires searching for augmenting paths multiple times.

\subsection{Conditions for Correct Solutions from Greedy Approach} \label{subsec:solutions}

To find conditions under which an algorithm that solves the SeqBAP can be used to solve the LexBAP with exactness guarantee, we explore the relationship between the solutions to the BAP, LexBAP, and SeqBAP.

\begin{proposition} \label{prop:sbap2bot}
	Consider a bipartite graph $\mathcal{G}_b$ and an MCM $\mathcal{M}$ of $\mathcal{G}_b$. If $\mathcal{M}$ is a solution to the SeqBAP as given in~(\ref{eq:SBAP}), then $\mathcal{M}$ is a solution to the BAP as given in~(\ref{eq:BAP}).
\end{proposition}
\begin{pf}
	Let $\mathcal{M}$ be a solution to~(\ref{eq:SBAP}). From~(\ref{eq:seqedgek}) and~(\ref{eq:M_k}), $\{i_1,j_1\}\in\mathcal{L}(\mathcal{M}^1)$, where $\mathcal{M}^1$ is a solution to~(\ref{eq:BAP}). By Lemma~\ref{lem:orderedges} in the appendix, $w(\{i_1,j_1\}) \geq w(e)$ for all edges $e\in\mathcal{M}$. Thus, $\mathcal{M}$ is also a solution to~(\ref{eq:BAP}) since $\{i_1,j_1\}\in\mathcal{L}(\mathcal{M}^1)$ and $\{i_1,j_1\}\in\mathcal{L}(\mathcal{M})$. $\hfill\qed$
\end{pf}

Combining the statements from Propositions~\ref{prop:priceBOT} and~\ref{prop:sbap2bot}, we obtain the following corollary. While Theorem~\ref{theorem:pos} provides a method to identify an edge with positive price of absence, Corollary~\ref{cor:priceLEX} motivates the need to identify edges with positive price of absence by showing their relevance in connection to the SeqBAP.

\begin{corollary} \label{cor:priceLEX}
	Consider a bipartite graph $\mathcal{G}_b=(\mathcal{V}_b,\mathcal{E}_b)$. If an edge $e_p\in\mathcal{E}_b$ has positive price of absence with respect to $\mathcal{G}_b$, then every MCM $\mathcal{M}$ that is a solution to the SeqBAP given in~(\ref{eq:SBAP}) contains $e_p$, i.e., $e_p\in\mathcal{M}$.
\end{corollary}

\begin{remark}
	Since Corollary~\ref{cor:priceLEX} also applies for each consecutive bipartite graph $\mathcal{G}^k$ in~(\ref{eq:G_k}), the following holds. Instead of selecting and removing one edge per iteration of $k$ in~(\ref{eq:seqedgek}), (\ref{eq:A_k}), (\ref{eq:B_k}), all edges with positive price of absence in that iteration can be selected and removed as a batch. Intuitively, edges are ``locked'' into the solution whenever they are found to have positive price of absence.
\end{remark}

Given Proposition~\ref{prop:sbap2bot}, we can go one step further; the following proposition states that all solutions to the LexBAP are also solutions to the SeqBAP.

\begin{proposition} \label{prop:lex2sbap}
	Consider a bipartite graph $\mathcal{G}_b$ and an MCM $\mathcal{M}$ of $\mathcal{G}_b$. If $\mathcal{M}$ is a solution to the LexBAP given in~(\ref{eq:LBAP}), then $\mathcal{M}$ is also a solution to the SeqBAP given in~(\ref{eq:SBAP}).
\end{proposition}
\begin{pf}
	Assume $\mathcal{M}=\{e_1,e_2,...,e_n\}$ is a solution to~(\ref{eq:LBAP}). Without loss of generality assume that $w(e_1) \geq w(e_2) \geq ... \geq w(e_n)$. Let $\mathcal{M}^1=\mathcal{M}$ and $\mathcal{M}^{k+1}=\mathcal{M}^k\backslash\{e_k\}$ for all $k\in\{1,2,...,n-1\}$. Then, it holds that for all $k\in\{1,2,...,n\}$, $e_k\in\mathcal{L}(\mathcal{M}^k)$. Since $\mathcal{M}$ is a solution to~(\ref{eq:LBAP}), for all $k\in\{1,2,...,n\}$, $\mathcal{M}^k \in \arg\!\min_{\mathcal{M}\in \mathcal{C}(\mathcal{G}^k)}\max_{\{i,j\}\in \mathcal{M}} w(\{i,j\})$, with $\mathcal{G}^k$ defined as in~(\ref{eq:G_k}). Then, it remains to show that $e_k\in\arg\!\max_{\{i,j\}\in \mathcal{L}(\mathcal{M}^k)}P(\mathcal{G}^k,\{i,j\})$. If $w(e_k) > w(e_{k+1})$, then $\mathcal{L}(\mathcal{M}^k)$ is a singleton, so $e_k$ is trivially the edge with largest price of absence in $\mathcal{G}^k$. If $\mathcal{L}(\mathcal{M}^k)$ is not a singleton, then it holds that $\mathcal{L}(\mathcal{M}^k)=\{e_k,e_{k+1},...,e_{k+\alpha}\}$, where either $k+\alpha=n$ or $w(e_{k+\alpha})>w(e_{k+\alpha+1})$.
	For the LexBAP solution, the choice of which $e\in\mathcal{L}(\mathcal{M}^k)$ is denoted as $e_k$ is arbitrary; edges $e_k,e_{k+1},...,e_{k+\alpha}$ can be rearranged in any order because their weights are equal. We choose $e_k= \arg\max_{e\in\{e_k,e_{k+1},...,e_{k+\alpha}\}} P(\mathcal{G}^k,e)$. Thus for all $k\in\{1,2,...,n\}$, $e_k\in\arg\max_{\{i,j\}\in \mathcal{L}(\mathcal{M}^k)}P(\mathcal{G}^k,\{i,j\})$ by construction as required in~(\ref{eq:seqedge}). $\hfill\qed$
\end{pf}

\begin{figure}[thpb]
	\centering
	\begin{subfigure}[b]{0.2\textwidth}
		\centering
		\begin{tikzpicture}[scale=0.5]
		\draw[, draw = black,fill opacity = .2] (0,0) circle (3.6);
		\draw[, draw = black,fill opacity = .2] (0,-1) circle (2.5);
		\draw[, draw = black,fill opacity = .2] (0,-1.5) circle (1.7);
		\node at (0,3.1) {$BAP$};
		\node at (0,0.75) {$SeqBAP$};
		\node at (0,-1) {$LexBAP$};
		\end{tikzpicture}
	\caption{}\label{fig:venn1}
	\end{subfigure}
	\hfill
	\begin{subfigure}[b]{0.2\textwidth}
		\centering
		\begin{tikzpicture}[scale=0.5]
			\draw[, draw = black,fill opacity = .2] (0,0) circle (3.6);
			\node at (0,-1) [circle,fill,inner sep=1.5pt]{};
			\node at (0,3.1) {$BAP$};
			\node at (0,-0.5) {$SeqBAP$, $LexBAP$};
		\end{tikzpicture}
	\caption{}\label{fig:venn2}
	\end{subfigure}
	\caption{Venn diagrams of the sets of the solutions to the BAP given in~(\ref{eq:BAP}), the SeqBAP given in~(\ref{eq:SBAP}), and the LexBAP given in~(\ref{eq:LBAP}). The diagram in (a) is for an arbitrary bipartite graph $\mathcal{G}_b$; this illustrates the general results provided by Propositions~\ref{prop:sbap2bot} and~\ref{prop:lex2sbap}. The diagram in (b) is for the case when there exist an MCM for which all edges have positive price of absence; this illustrates Corollary~\ref{cor:unique}.}\label{fig:venn}
\end{figure}
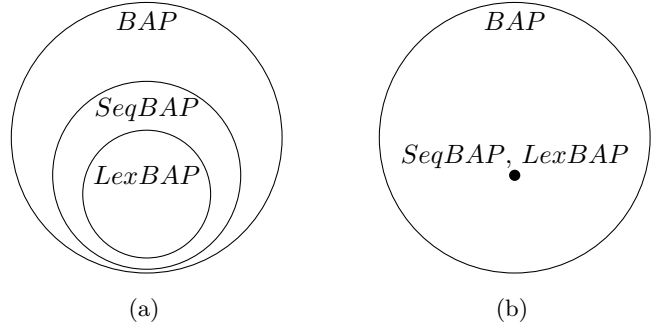

Fig.~\ref{fig:venn}(a) illustrates Propositions~\ref{prop:sbap2bot} and~\ref{prop:lex2sbap}. The following corollary is a special case of Proposition~\ref{prop:lex2sbap}, where the solution set of the SeqBAP is a singleton.

\begin{corollary} \label{cor:unique}
	Consider a bipartite graph $\mathcal{G}_b$ and an MCM $\mathcal{M}$ of $\mathcal{G}_b$. If $\mathcal{M}$ is the unique solution to the SeqBAP given in~(\ref{eq:SBAP}), then $\mathcal{M}$ is also the unique solution to the LexBAP given in~(\ref{eq:LBAP}).
\end{corollary}

The converse of Corollary~\ref{cor:unique} is not true. We provide the following counterexample with two agents $i_1, i_2$ and two tasks $j_1, j_2$, and edges with weights $w(\{i_1,j_1\})=2$, $w(\{i_1,j_2\})=2$, $w(\{i_2,j_1\})=1$, $w(\{i_2,j_2\})=2$. In this case, the LexBAP has a unique solution, but the SeqBAP does not. Fig.~\ref{fig:venn}(b) illustrates Corollary~\ref{cor:unique}.


The following proposition provides conditions for existence of a unique solution to the SeqBAP.

\begin{proposition} \label{prop:pos}
	Consider a bipartite graph $\mathcal{G}_b$ and an MCM $\mathcal{M}$ of $\mathcal{G}_b$. Let $\mathcal{M}$ be the solution to the SeqBAP given in~(\ref{eq:SBAP}). All sequentially selected edges in~(\ref{eq:seqedgek}) have positive price of absence in their respective graphs $\mathcal{G}^k$ defined in~(\ref{eq:G_k}) if and only if $\mathcal{M}$ is a unique solution to~(\ref{eq:SBAP}).
\end{proposition}
\begin{pf}
	Assume all sequentially selected edges in~(\ref{eq:seqedgek}) have positive price of absence in their respective graphs $\mathcal{G}^k$. By Corollary~\ref{cor:priceLEX}, every solution to~(\ref{eq:SBAP}) must contain this set of edges. Thus, the solution is unique.
	
	Assume there exists an edge $e_q$ selected in~(\ref{eq:seqedgek}) at iteration $q$ that does not have positive price of absence in $\mathcal{G}^q$. Then there exists another bottleneck edge of $\mathcal{G}^q$ that can be selected in lieu of $e_q$ at iteration $q$. Thus, the solution to~(\ref{eq:SBAP}) is not unique. $\hfill\qed$
\end{pf}

Proposition~\ref{prop:pos} is stronger than the result in~\cite{tony}, which only considers the sufficiency but the not necessity of all edges having positive price for uniqueness of the SeqBAP solution. The relationship between the solution sets of the BAP, SeqBAP and LexBAP allows us to re-derive the following result from~\cite{tony}. In particular, by combining Corollary~\ref{cor:unique} and Proposition~\ref{prop:pos}, we have a sufficient condition for a solution to the LexBAP being unique.

\begin{corollary} \label{cor:posimpliesunique}
	Consider a bipartite graph $\mathcal{G}_b$ and an MCM $\mathcal{M}$ of $\mathcal{G}_b$. Let $\mathcal{M}$ be the solution to the SeqBAP given in~(\ref{eq:SBAP}). If all sequentially selected edges in~(\ref{eq:seqedgek}) have positive price of absence in their respective graphs $\mathcal{G}^k$ defined in~(\ref{eq:G_k}), then $\mathcal{M}$ is a unique solution to the LexBAP given in~(\ref{eq:LBAP}).
\end{corollary}

We in turn note that a sufficient condition for all sequentially selected edges in a solution to the SeqBAP having positive price of absence is for all weights in the bipartite graph to be distinct. Therefore, if the weights of all edges in $\mathcal{G}_b$ are distinct, then the SeqBAP has a unique solution and this solution is also the unique solution of the LexBAP. In realistic applicatons, the weights can often be considered to belong to a non-empty interval of real numbers. This occurs for instance in applications where the weights consist of distances between agents and tasks. In such applications, the situation where the weights are non-distinct has zero measure.


\section{A Distributable Greedy Algorithm for Solving the LexBAP} \label{sec:alg}

We present a method to solve the LexBAP that exploits the structure analysed in Section~\ref{sec:structs}. In particular, we introduce Algorithm~\ref{alg:lex} that solves the SeqBAP. The algorithm provides certificates for when this solution to the SeqBAP is also a solution to the LexBAP and it can be implemented with distributed computation across agents. The following assumptions model a distributed setting, where the information available to each individual agent may be limited.

\begin{assumption} \label{as:distedge}
	Assume an agent $i\in\mathcal{A}_b$ has access to the set of incident edges $\mathcal{E}_i := \{\{i,j\}\in \mathcal{E}_b|j\in\mathcal{B}_b\}$ and access to the weight of each edge in $\mathcal{E}_i$, i.e., the set $C_i:=\{w(\{i,j\})|\{i,j\}\in\mathcal{E}_i\}$.
\end{assumption}

Note, $\mathcal{E}_b=\bigcup_{i\in\mathcal{A}_b} \mathcal{E}_i$ and $\mathcal{E}_{v}\cap \mathcal{E}_{v'}=\emptyset$ for $v,v'\in\mathcal{A}_b$, $v\neq v'$.


\begin{assumption} \label{as:distcom}
	Let communication between agents be modelled by a time invariant, undirected and connected graph $\mathcal{G}_C=(\mathcal{A}_b,\mathcal{E}_C)$ with vertex set $\mathcal{A}_b$, edge set $\mathcal{E}_C$ and diameter $D$. Assume all agents $i\in\mathcal{A}_b$ communicate synchronously according to a global clock.
\end{assumption}

Assumption~\ref{as:distcom} describes the communication between agents. Synchronous communication refers to communication where all agents $i\in\mathcal{A}_b$ share a global clock and at each time step of the clock, agents exchange information with all their neighbours $i'\in N(\mathcal{G}_C,i)$. Thus, there is a delay for information to reach a non-neighbouring agent as that information is relayed through a chain of agents. This delay is at most $D$ time steps of the global clock, corresponding to the diameter of $\mathcal{G}_C$. Algorithm~\ref{alg:lex} requires Assumption~\ref{as:distcom}, for example, when agents need to reach a consensus on the global edge with largest weight in $\mathcal{M}$ in line~\ref{line:max}. In this case, agents must make use of communication to reach a global consensus on the edge with largest weight as each agent is only initially aware of its own local edge with largest weight.


Algorithm~\ref{alg:lex} can be initialised with any arbitrary MCM $\mathcal{M}_0$. Algorithm~\ref{alg:lex} uses the function $\texttt{\textsc{AugPath}}(\cdot)$ introduced in~\cite{distBAP} to solve the SeqBAP. Consider a bipartite graph $\mathcal{G}_b=(\mathcal{V}_b,\mathcal{E}_b)$, an MCM $\mathcal{M}$, and an edge $e\in\mathcal{M}$. If an augmenting path exists relative to $\mathcal{M}\backslash \{e\}$ within the set $\mathcal{E}\subseteq\mathcal{E}_b$, then the function $\texttt{\textsc{AugPath}}(\mathcal{M}\backslash \{e\},\mathcal{E})$ returns an MCM $\mathcal{M}'\subseteq \mathcal{E}$. If an augmenting path does not exist, then it returns $\mathcal{M}\backslash \{e\}$. In~\cite{distBAP}, it is shown that $\texttt{\textsc{AugPath}}(\cdot)$ can be implemented with the distributed setting given in Assumptions~\ref{as:distedge} and~\ref{as:distcom}. The function $\texttt{\textsc{AugPath}}(\cdot)$ in Algorithm~\ref{alg:lex} requires the following assumption on the input graph.

\begin{assumption}\label{as:AgB}
	The graph $\mathcal{G}_b=(\mathcal{A}_b\cup\mathcal{B}_b,\mathcal{E}_b)$ has cardinalities of agent and task sets satisfiying $|\mathcal{A}_b|\geq |\mathcal{B}_b| = n$, where $n$ is the cardinality of an MCM of $\mathcal{G}_b$, i.e., the number of agents is greater than or equal to the number of tasks.
\end{assumption}

\makeatletter
\algrenewcommand\ALG@beginalgorithmic{\footnotesize}
\makeatother

\begin{algorithm}[thpb]
	\caption{Algorithm for solving the SeqBAP.} \label{alg:lex}
	Input: Graph $\mathcal{G}_b=(\mathcal{V}_b,\mathcal{E}_b)$ and an MCM $\mathcal{M}_0$.\\
	Output: An MCM $\mathcal{M}_b$ of $\mathcal{G}_b$ that is a solution to~(\ref{eq:SBAP}) and a flag $\mathcal{U}$ for it being an exact solution to~(\ref{eq:LBAP}).
	\begin{algorithmic}[1]
		\State $\bar{\mathcal{V}}\gets\mathcal{V}_b$
		\State $\bar{\mathcal{E}}\gets\mathcal{E}_b$
		\State $\bar{\mathcal{M}}\gets \mathcal{M}_0$
		\State $\mathcal{M}_b\gets \emptyset$
		\State $\mathcal{U}\gets \texttt{True}$
		\While {$|\bar{\mathcal{M}}|>0$} \label{line:while}
		\State $\bar{\mathcal{G}}\gets (\bar{\mathcal{V}},\bar{\mathcal{E}})$\Comment{{\scriptsize Current graph}} \label{line:start}
		\State $\bar{e}\gets e\in\mathcal{L}(\bar{\mathcal{M}})$ \label{line:max} \Comment{{\scriptsize Find edge with largest weight in $\bar{\mathcal{M}}$}}
		\State $\bar{\mathcal{E}}\gets \bar{\mathcal{M}}\cup\psi_S(\bar{\mathcal{G}},\bar{\mathcal{M}})$ \Comment{{\scriptsize Shrink $\bar{\mathcal{E}}$}}
		\State $\mathcal{M}_\nu \gets \texttt{\textsc{AugPath}}(\bar{\mathcal{M}}\backslash \{\bar{e}\},\bar{\mathcal{E}}\backslash \{\bar{e}\})$\label{line:aug1}
		\If {$\mathcal{M}_\nu \neq \bar{\mathcal{M}}$} \Comment{{\scriptsize $\bar{e}$ is not a critical bottleneck edge}}
		\State $\bar{\mathcal{M}}\gets \mathcal{M}_\nu$
		\Else \Comment{{\scriptsize $\bar{e}$ is a critical bottleneck edge}} \label{line:end}
		\State $E'\gets \emptyset$
		\For {$e'\in\bar{\mathcal{M}}$}
		\State $\mathcal{M}_\nu \gets\texttt{\textsc{AugPath}}(\bar{\mathcal{M}}\backslash \{e'\},\psi(\bar{\mathcal{G}},\bar{\mathcal{M}})\backslash \{e'\})$\label{line:pos}
		\If {$\mathcal{M}_\nu= \bar{\mathcal{M}}$} \Comment{{\scriptsize $P(\bar{\mathcal{G}},e')>0$}}
		\State $E'\gets E'\cup \{e'\}$ \Comment{{\scriptsize See Corollary~\ref{cor:E}}}
		\EndIf
		\EndFor
		\If {$E'\cap \mathcal{L}(\bar{\mathcal{M}})=\emptyset$} \Comment{{\scriptsize See Remark~\ref{rem:nopos}}} \label{line:zerorob}
		\State $\mathcal{U}\gets \texttt{False}$ \Comment{{\scriptsize See Proposition~\ref{prop:pos}}}
		\State $E'\gets E'\cup\{\bar{e}\}$ \label{line:atleastone}
		\EndIf
		\State $\mathcal{V}'\gets \{v\in\mathcal{V}_b|\{v,v'\}\in E'\}$
		\State $\bar{\mathcal{V}}\gets \bar{\mathcal{V}} \backslash \mathcal{V}'$\Comment{{\scriptsize Remove vertices incident to edges in $E'$}} \label{line:vertex}
		\State $\bar{\mathcal{E}}\gets \{\{i,j\}\in\bar{\mathcal{E}}|i\in\bar{\mathcal{V}} \text{ and }j\in\bar{\mathcal{V}}\}$ \Comment{{\scriptsize Shrink $\bar{\mathcal{E}}$}} 
		\State $\bar{\mathcal{M}} \gets \bar{\mathcal{M}} \backslash E'$ \label{line:prune}
		\State $\mathcal{M}_b\gets \mathcal{M}_b\cup E'$
		\EndIf
		\EndWhile
		\State \Return $\mathcal{M}_b$, $\mathcal{U}$
	\end{algorithmic}
\end{algorithm}

The algorithm runs by systematically removing elements from a graph that is initialised with $\mathcal{G}_b$. The following steps are repeated in each iteration of the while-loop beginning in Line~\ref{line:while}. First, one edge from the current graph is tested to see if it is a critical bottleneck edge. Testing to see if an edge is a critical bottleneck edge involves an augmenting path search. If a critical bottleneck edge is found, all edges with positive price of absence in the current MCM are identified. This can be implemented by again searching for augmenting paths, which can be carried out with $\texttt{\textsc{AugPath}}(\cdot)$. All edges with positive price of absence and all edges adjacent to these edges are removed from the current graph. If none of the bottleneck edges of the current graph have positive price of absence, an arbitrary bottleneck edge and the edges adjacent to it are removed in accordance with~(\ref{eq:seqedgek}) and (\ref{eq:A_1}-\ref{eq:B_k}), see Remark~\ref{rem:nopos}. This ensures the graph reduces by at least one edge each time a critical bottleneck edge is found. After making these changes, the next iteration of the while-loop is commenced on the reduced graph.

\begin{remark} \label{rem:nopos}
	If none of the bottleneck edges have positive price of absence, i.e., $E'\cap\mathcal{L}(\mathcal{M})=\emptyset$, where $E'$ is given in Corollary~\ref{cor:E} and $\mathcal{L}(\mathcal{M})$ is defined in~(\ref{eq:mathl}), then any arbitrary bottleneck edge satisfies~(\ref{eq:seqedgek}).
\end{remark}

A key observation is that the graph always reduces in each iteration by removal of edges that are not found to be critical bottleneck edges, by removal of edges with positive price of absence together with edges adjacent to them, or by removal of an edge that is found to be a critical bottleneck edge together with edges adjacent to it. The total number of iterations of the while-loop beginning in Line~\ref{line:while} is upper bounded by $|\mathcal{E}_b|$, and depends on how quickly the pool of candidate critical bottleneck edges $\bar{\mathcal{E}}$ shrinks and how many edges with positive prices of absence $|E'|$ are found in each iteration.

\begin{theorem}
	Algorithm~\ref{alg:lex} produces a solution to the SeqBAP and can be implemented in the distributed setting given by Assumptions~\ref{as:distedge},~\ref{as:distcom} and~\ref{as:AgB}.
\end{theorem}
\begin{pf}
	We first prove that Algorithm~\ref{alg:lex} is amenable to implementation in the distributed setting. To this end, we observe that Algorithm~\ref{alg:lex} only requires three procedures, i.e., the procedures of edge removal, finding the edge with largest weight and searching for an augmenting path. All three procedures have been shown to be distributable in~\cite{distBAP} and can therefore be implemented in the distributed setting.
	
	Next, we apply the results from Section~\ref{sec:structs} to prove convergence to a SeqBAP solution. A bottleneck assignment is found in Lines~\ref{line:start} to~\ref{line:end}, the proof of which can be found in~\cite{distBAP}. By Corollary~\ref{cor:E}, we can identify multiple edges in a SeqBAP solution by checking their positivity of price of absence. By Theorem~\ref{theorem:pos}, this involves augmenting path searches as indicated by Line~\ref{line:pos}. The graph is pruned in Lines~\ref{line:zerorob} to~\ref{line:prune} according to how many edges have positive price of absence, i.e., the number of edges in~\eqref{eq:seqedge} that have been selected so far. Line~\ref{line:atleastone} ensures that at least one edge is selected in each loop of Algorithm~\ref{alg:lex}. The process is repeated until all $n$ edges in~\eqref{eq:seqedge} have been selected. $\hfill\qed$
\end{pf}

Apart from an MCM, the algorithm returns a flag when this MCM is a unique solution to the LexBAP in accordance with Corollary~\ref{cor:posimpliesunique}. If the guard in Line~\ref{line:zerorob} is false for all iterations, then the solution is unique as described in Corollary~\ref{cor:unique}. If the flag is not returned as true, then the SeqBAP has multiple solutions and the produced MCM may not be a solution to the LexBAP. 




\begin{proposition} \label{prop:complex}
	Assume the number of agents and the number tasks in $\mathcal{G}_b$ are both equal to $n$. The worst-case complexity of Algorithm~\ref{alg:lex} is $\mathcal{O}(n^3D)$, where $D$ is the diameter of the communication graph.
\end{proposition}
\begin{pf}
	At most $n^2$ edges are tested as candidate critical bottleneck edges in Line~\ref{line:aug1} of Algorithm~\ref{alg:lex}. Finding and testing an edge involves a distributed max-consensus and an augmenting path search. In this setting, these procedures have orders $\mathcal{O}(D)$ and $\mathcal{O}(nD)$ respectively, as shown in~\cite{distBAP}. The function $\texttt{\textsc{AugPath}}(\cdot)$ has complexity $\mathcal{O}(nD)$ because it exploits the fact that a matching of size $n-1$, and the most recently removed edge $\bar{e}$ are both known inputs. In the worst-case, every edge in $\mathcal{E}_b$ is tested in this way. Therefore, the worst-case complexity of applying these two procedures to every edge is $\mathcal{O}(n^3D)$. No edge is tested for being a critical bottleneck edge more than once, and at most $n$ critical bottleneck edges must be found.
	
	Each time a candidate proves to be a critical bottleneck edge, the test in Line~\ref{line:pos} is carried out to identify edges with positive price of absence. Testing one edge for positive price according to Theorem~\ref{theorem:pos} involves an augmenting path search, which as mentioned, has complexity $\mathcal{O}(nD)$. By the contrapositive of Proposition~\ref{prop:priceBOT}, only edges in the current MCM $\mathcal{M}$ are candidates that need to be tested and $\mathcal{M}$ has at most $n$ edges. The complexity of testing all edges in an MCM that has a maximum cardinality of $n$ is $\mathcal{O}(n^2D)$.
	
	Testing all edges for being critical bottleneck edges has complexity $\mathcal{O}(n^3D)$. Testing edges for positive price of absence has complexity $\mathcal{O}(n^2D)$ per MCM, but it is done at most $n$ times corresponding to the maximum number of critical bottleneck edges, so in the worst-case it is also $\mathcal{O}(n^3D)$. Thus, the complexity of Algorithm~\ref{alg:lex} is $\mathcal{O}(n^3D)$. $\hfill\qed$
\end{pf}


As a comparison, the centralised algorithm to solve the LexBAP presented in \cite{assignprob} has complexity $\mathcal{O}(n^4)$ and involves $n$ iterations of solving both a BAP with a complexity of $\mathcal{O}(n^{2.5})$~\cite{assignprob} and a Linear Sum Assignment Problem (LSAP)~\cite{hungarian,LAP1} with a complexity of $\mathcal{O}(n^3)$.

With $D=1$, Algorithm~\ref{alg:lex} has a worst-case complexity of $\mathcal{O}(n^3)$; $D=1$ is a special case corresponding to a centralised algorithm as agents communicate with to all other agents. The more complex LSAPs are bypassed by applying the augmenting path searches in this greedy approach, yet under the conditions described in Section~\ref{subsec:solutions} the solution found using either algorithm is identical.

A naive greedy approach for solving the LexBAP that finds $n$ bottleneck edges from scratch, where each subsequent bottleneck edge is found without exploiting knowledge of previous bottleneck assignments, has complexity $\mathcal{O}(n^{3.5})$. This type of naive greedy approach using an ``off-the-shelf'' BAP algorithm in sequence was first proposed in~\cite{tony,seq1} . On a related note, the complexity of a SeqBAP algorithm that additionally returns the explicit prices of absence of edges is also $\mathcal{O}(n^{3.5})$ and requires solving $n$ iterations of two BAPs. While Algorithm~\ref{alg:lex} relies on identifying edges with positive price of absence, the explicit value of the price of absence is never computed. In applications that utilise the value of the price of absence, e.g., to quantify the robustness of an assignment in~\cite{tony}, additional computation is required.

\section{Case Study} \label{sec:case}

Consider $n$ agents represented by points in $\mathbb{R}^2$. Each agent must move from its initial position to one of $n$ goal positions. The assignment of agents to goal positions is done by solving the LexBAP, where the weights are given by the agent-goal distances. All coordinates are generated from a uniform distribution between values of 0 and 100 normalised distance units. Since the weights are almost surely distinct, the solution to the SeqBAP is the solution to the LexBAP and it is unique.

An example application for this case study would be a ride-sharing service, where the agents are cars and the goal positions are pick-up locations. Another example is a drone delivery service or autonomous forklifts in a warehouse, where the agents are drones or forklifts and the goal positions are item pick-up locations. Assuming tasks are executed in parallel, all these applications correspond to a BAP as the longest pick-up time needs to be minimised. In contrast, minimising the sum of costs may result in some agents completing their tasks very quickly at the expense of other agents completing theirs very slowly. Furthermore from Fig.~\ref{fig:venn}, solving the BAP allows the possibility of further minimising the sequential bottleneck edges, i.e., solving the LexBAP. The benefit of this is that the LexBAP itself has further desirable properties, e.g., an inherent property that provides collision avoidance guarantees as proven in~\cite{tony}.

\begin{figure}[thpb]
	\centering
	\includegraphics[scale=0.26]{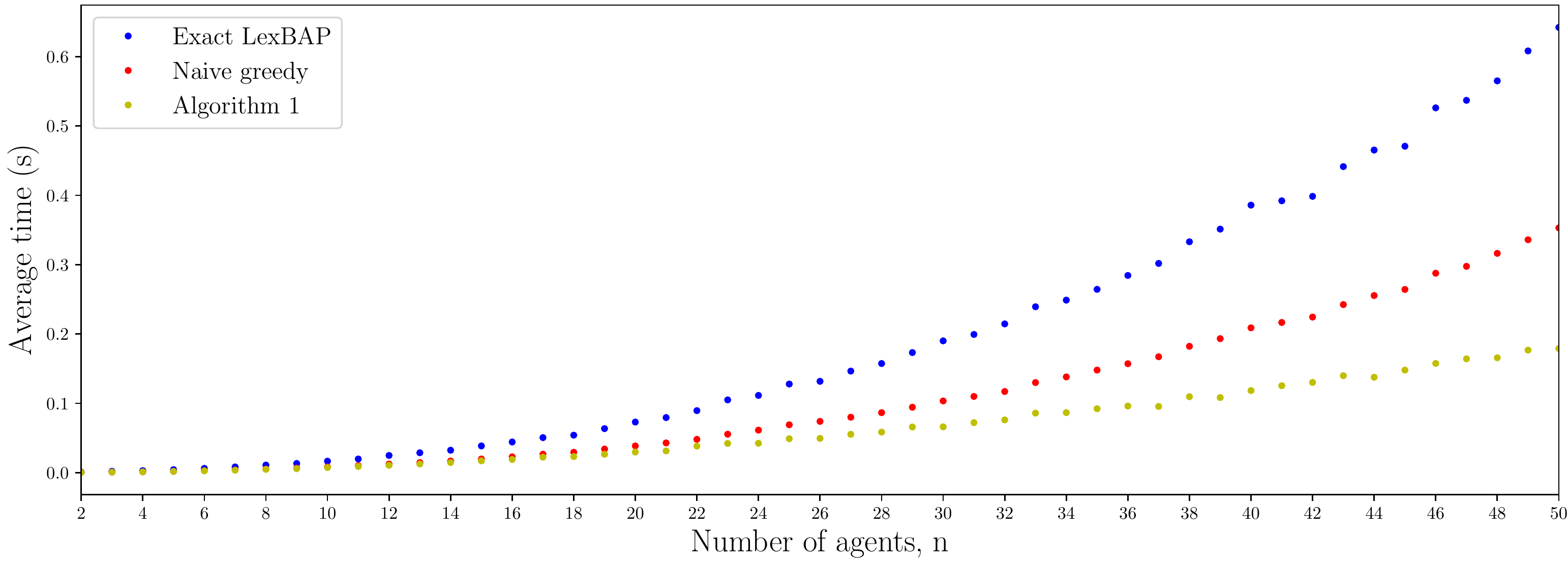}
	\caption{Comparison of the runtime of an exact LexBAP algorithm, a naive greedy approach, and Algorithm~\ref{alg:lex}.} \label{fig:time}
\end{figure}

Fig.~\ref{fig:time} shows the performance of a centralised version of Algorithm~\ref{alg:lex} benchmarked against an algorithm for solving general LexBAPs and the naive greedy approach for solving the LexBAP. The LSAP component of the exact LexBAP method is solved using the Hungarian Algorithm from~\cite{hungarian} while the BAP component in both the exact LexBAP and the naive greedy approach is solved using an off-the-shelf threshold method from~\cite{assignprob}. The average time for each value of $n$ is taken over 100 realisations of agent and task positions. These results are obtained with an Intel i5-6600 CPU at 3.30GHz. For all algorithms tested, Fig.~\ref{fig:time} shows the average time increases as the number of agents and tasks in the problem increases, i.e., as $n$ increases.

From Proposition~\ref{prop:complex}, the theoretical worst-case complexity of Algorithm~\ref{alg:lex} is $\mathcal{O}(n^3D)$, which is lower than the theorectical worst-case complexities of the exact LexBAP and naive greedy approach that are order $\mathcal{O}(n^4D)$ and $\mathcal{O}(n^{3.5}D)$ respectively. On the other hand, Fig.~\ref{fig:time} shows that the empirical complexity of Algorithm~\ref{alg:lex} is also lower than the empirical complexities the exact LexBAP and naive greedy approach for this case study.

\begin{figure}[thpb]
	\centering
	\includegraphics[scale=0.4]{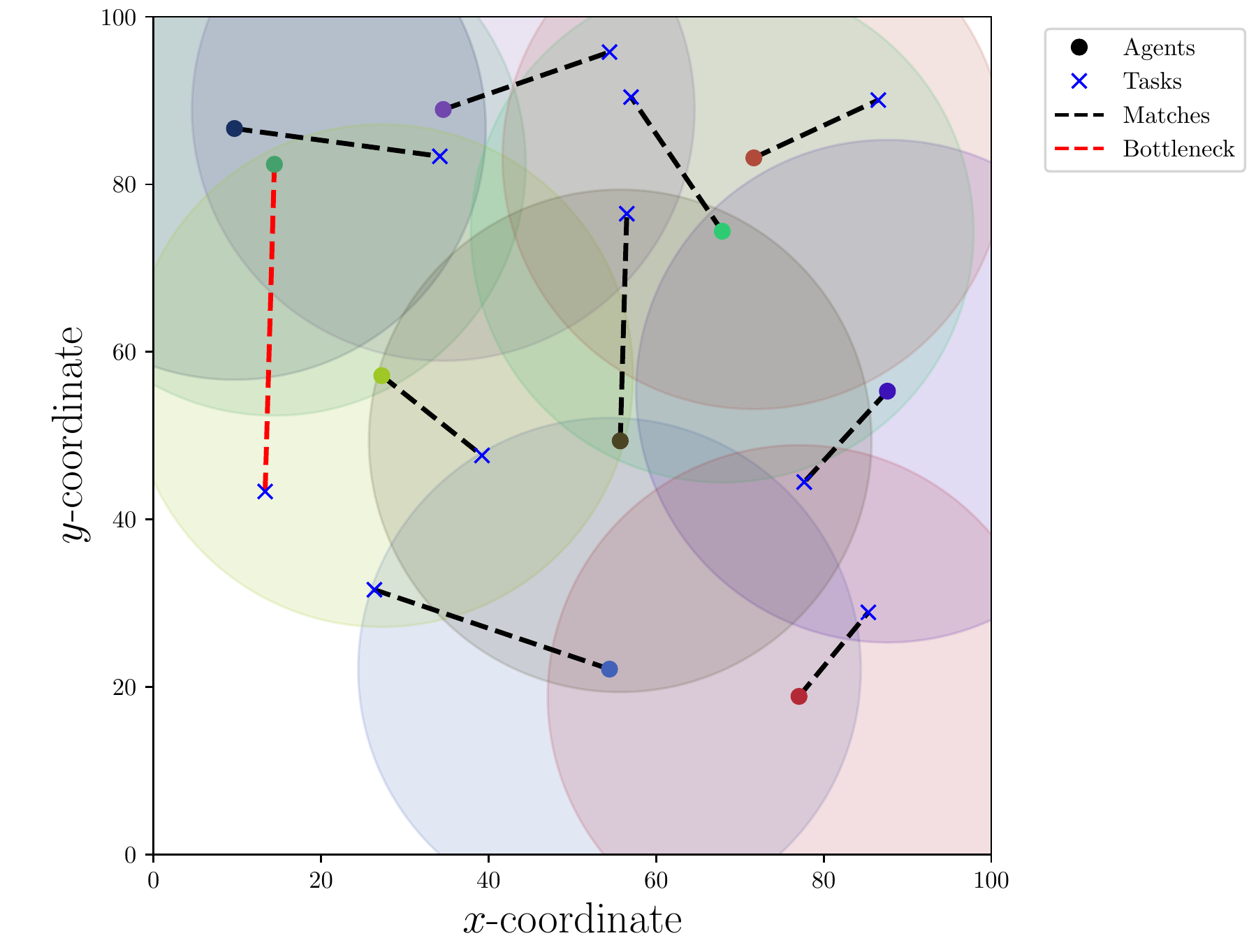}
	\caption{Demonstration of Algorithm~\ref{alg:lex} in a distributed setting. The shaded circles show the communication range of each agent. Agents are only able to communicate with other agents within their circles.} \label{fig:diam}
\end{figure}

The following demonstrates that Algorithm~\ref{alg:lex} can be implemented under the distributed setting given by Assumptions~\ref{as:distedge} and~\ref{as:distcom}. Fig.~\ref{fig:diam} shows one realisation of agent and task positions with $n=10$. First, Algorithm~\ref{alg:lex} is applied for the case where any agent can communicate with all other agents directly in one time step of the global clock, i.e., the communication graph has diameter $D = 1$. Additionally, we consider the case where agents only communicate with other agents that are located within a radius of $30$ units as illustrated by the shaded areas in Fig.~\ref{fig:diam}. This results in a communication graph with diameter $D = 5$. For both cases, Alogorithm~\ref{alg:lex} returns the exact the solution of the LexBAP. For the case with $D = 1$, it takes 111 time steps for Algorithm~\ref{alg:lex} to obtain the solution. For the case with $D = 5$, it takes 555 time steps.


\section{Conclusion} \label{sec:conclusion}


We presented an approach to find an MCM of a bipartite graph that is the solution to the LexBAP by employing a method that solves a series of BAPs, where the edges in the bipartite graph are removed in each iteration. For each of these BAPs, we showed that if an edge has a positive price of absence, then that edge is guaranteed to be an element of the LexBAP solution, and there may be multiple such edges each time a BAP is solved. We called this greedy reformulation of the LexBAP the SeqBAP. We considered the similarities in structure of a critical bottleneck edge to an edge with positive price of absence and used this to derive a method to identify edges with positive price of absence that involves a search for augmenting paths. This enables solving the SeqBAP efficiently. We derived the relationship between the BAP, the SeqBAP, and the LexBAP by comparing their solution sets. In particular, we showed that the solutions to the LexBAP are a subset of the solutions of the SeqBAP. We provided conditions for when the SeqBAP has a unique solution, which implies that this solution also uniquely solves the LexBAP. Futhermore, we showed that all edges of the graph having distinct weight values is a sufficient condition for this uniqueness. We combined these results into a proposed algorithm that provides a greedy solution to the LexBAP and a certificate for when this solution is exact. The algorithm has a complexity of $\mathcal{O}(n^3)$, which is lower than methods for solving the LexBAP that have a complexity of $\mathcal{O}(n^4)$. Moreover, the proposed algorithm can be implemented with computation that is distributed across agents.


\bibliographystyle{plain}
\bibliography{biblio}

\appendix
\section{Order of Edges Chosen by SeqBAP}

\begin{lemma} \label{lem:orderedges}
	The edges in~(\ref{eq:seqedge}) have weights such that $w(\{i_k,j_k\})\geq w(\{i_{k+1},j_{k+1}\})$, for all $k\in\{1,2,...,n-1\}$.
\end{lemma}
\begin{pf}
	For all $k\in\{1,2,...,n-1\}$, $w(\{i_k,j_k\})\geq w(e_k)$ for any $e_k\in\mathcal{M}^k$ because from~(\ref{eq:seqedge}), $\{i_k,j_k\}$ is selected from the set $\mathcal{L}(\mathcal{M}^k)$. The matching $\mathcal{M}^k\backslash \{\{i_k,j_k\}\}$ is an MCM of $\mathcal{G}^{k+1}$ but not necessarily the solution to the BAP for the graph $\mathcal{G}^{k+1}$. On the other hand, by~(\ref{eq:M_k}), $\mathcal{M}^{k+1}$ is a solution to the BAP for the graph $\mathcal{G}^{k+1}$. Thus, $w(e_{k+1})\leq w(e)$, where $e_{k+1}\in\mathcal{L}(\mathcal{M}^{k+1})$ and $e\in\mathcal{L}(\mathcal{M}^k\backslash \{\{i_k,j_k\}\})$. Therefore, $w(\{i_k,j_k\})\geq w(e)\geq w(e_{k+1})=w(\{i_{k+1},j_{k+1}\})$.$\hfill\qed$
\end{pf}


\end{document}